\documentclass{amsart}
\usepackage{amssymb} 
\newtheorem{thm}{Theorem}[section]
\newtheorem{cor}[thm]{Corollary}

\theoremstyle{definition}

\theoremstyle{remark}

\numberwithin{equation}{section}

\begin{document}
\title[Groups covers and blocking sets]{Groups with maximal irredundant covers and minimal blocking sets}%
\author{Alireza Abdollahi}%
\address{Department of Mathematics, University of Isfahan, Isfahan 81746-73441, Iran; and  School of Mathematics, Institute for Research in Fundamental Sciences (IPM), P.O.Box: 19395-5746, Tehran, Iran. e-mail: {\tt
abdollahi@member.ams.org} \;\; {\tt
a.abdollahi@math.ui.ac.ir} }

\thanks{}%

\begin{abstract}
Let $n$ be a positive integer. Denote by $\mathrm{PG}(n,q)$  the
$n$-dimensional projective space over the finite field
$\mathbb{F}_q$ of order $q$. A blocking set  in $\mathrm{PG}(n,q)$ is a set
of points that has non-empty intersection with every hyperplane of
$\mathrm{PG}(n,q)$. A blocking set is called minimal if none of its proper
subsets are blocking sets.
In this note we prove that if  $\mathrm{PG}(n_i,q)$ contains a minimal blocking set of size $k_i$ for
$i\in\{1,2\}$, then $\mathrm{PG}(n_1+n_2+1,q)$ contains a minimal blocking
set of size $k_1+k_2-1$. This  result is  proved by a result on groups with maximal irredundant covers.
\end{abstract}
\maketitle
\section{\bf Introduction and Results}
Let $G$ be a group. A set $\mathcal{S}$ of proper subgroups of $G$
is called a cover for $G$ whenever $G=\bigcup_{H\in \mathcal{S}}
H$. The cover $\mathcal{S}$ is called irredundant if no proper
sub-collection of $\mathcal{S}$ is a cover for $G$. The cover
$\mathcal{S}$ is called an $n$-cover if $|\mathcal{S}|=n$. The
cover $\mathcal{S}$ is called maximal if each member of
$\mathcal{S}$ is a maximal subgroup of the group $G$. The cover
$\mathcal{S}$ is called core-free if the core of $D=\bigcap_{H\in
\mathcal{S}}H$ in $G$ is trivial, i.e. $D_G=\bigcap_{g\in G}
g^{-1}Dg$ is the trivial subgroup of $G$. The cover $\mathcal{S}$
for a group is called a $\mathfrak{C}_n$-cover if $\mathcal{S}$
is a maximal irredundant core-free $n$-cover.

Let $n$ be a positive integer. Denote by $\mathrm{PG}(n,q)$  the
$n$-dimensional projective space over the finite field
$\mathbb{F}_q$ of order $q$. A blocking set  in $\mathrm{PG}(n,q)$ is a set
of points that has non-empty intersection with every hyperplane of
$\mathrm{PG}(n,q)$. A blocking set is called minimal if none of its proper
subsets are blocking sets.

There is a well-known relationship between minimal blocking sets
in $\mathrm{PG}(n,q)$ and irredundant covers of the abelian group
$\mathbb{F}_q^{n+1}$, the direct product of $n+1$ copies of
$\mathbb{F}_q$ (see e.g., \cite{AAM}). In fact a minimal blocking
set of size $k$ exits in $\mathrm{PG}(n,q)$ if and only if an irredundant
$k$-cover for $\mathbb{F}_q^{n+1}$ whose all members are
hyperplane exists (see e.g., Proposition 2.2 of \cite{AAM}).

The main results of this note are the following. Theorem \ref{thm1} is inspired by the
proof of Lemma 3.1 of \cite{T}.
\begin{thm}\label{thm1}
Let $\mathcal{C}_1=\{M_1,\dots,M_m\}$ and
$\mathcal{C}_2=\{N_1,\dots,N_n\}$ be irredundant $m$- and
$n$-covers for two groups $G_1$ and $G_2$, respectively. Suppose
that $M_1$ and $N_1$ are  normal subgroups of $G_1$ and $G_2$,
respectively; and assume that $a$ and $b$ are  elements of
$G_1\backslash M_1$ and $G_2\backslash N_1$, respectively, such
that $M_1a$ and $N_1b$ have the same order $p$ for some prime
number $p$. Then
$$\mathcal{C}=\{(M_1\times N_1)\langle (a,b)\rangle, M_i\times G_2,G_1\times
N_j \;|\; i=2,\dots,m; \; j=2,\dots,n\}$$ is an irredundant
$(n+m-1)$-cover for $G_1\times G_2$ with intersection $D_1\times
D_2$, where $D_1=\cap_{i=1}^m M_i$, $D_2=\cap_{j=1}^n N_i$ and $(M_1\times N_1)\langle (a,b)\rangle$ is the subgroup of $G_1\times G_2$ generated by
 $M_1\times N_1$ and the element $(a,b)$. In
particular, if both $\mathcal{C}_1$ and $\mathcal{C}_2$ are
maximal, then  $\mathcal{C}$ is a maximal cover for $G_1\times
G_2$.
\end{thm}
\begin{cor}\label{cor}
Let $\mathrm{PG}(n_i,q)$ contain a minimal blocking set of size $k_i$ for
$i\in\{1,2\}$. Then $\mathrm{PG}(n_1+n_2+1,q)$ contains a minimal blocking
set of size $k_1+k_2-1$.
\end{cor}
\section{\bf Proofs}
 \noindent{\bf Proof of Theorem \ref{thm1}.}
It is clear that $\mathcal{C}$ is  an $(m+n-1)$-cover for
$G_1\times G_2$. We now prove that $\mathcal{C}$ is irredundant.
 First, $(M_1\times N_1)\langle (a,b)\rangle$ is an irredundant member of $\mathcal{C}$: for $\mathcal{C}_1$ and $\mathcal{C}_2$ are both irredundant,
there exist  elements $x\in M_1\backslash \bigcup_{i=2}^m M_i$ and
$y\in N_1\backslash \bigcup_{j=2}^n N_i$. It follows that the
element $(x,y)$, among the members of $\mathcal{C}$, belongs only
to $(M_1\times N_1)\langle (a,b)\rangle$ and so $(M_1\times
N_1)\langle (a,b)\rangle$ is an irredundant member. Now we show that  each $M_i\times G_2$ is an
irredundant member for the cover (for each $i=2,\dots, m$). Without
loss of generality  and for convenience, we prove only the latter
statement for $i=2$, the other cases are similar. Let $x_0\in
M_2\backslash \bigcup_{\overset{i=1}{i\neq 2}}^m M_i$ and $y_0\in
N_1\backslash \bigcup_{j=2}^n N_j$. Then $(x_0,y_0)$, among the
members of $\mathcal{C}$, could only possibly belong  to $(M_1\times N_1)\langle
(a,b)\rangle$. If it is possible, then  $x_0 \in M_1a^\ell$ and $y_0\in
N_1b^\ell$ for some integer $\ell$. Since $x_0\not\in M_1$,
$a^\ell\not\in M_1$ and so $p$ does not divide $\ell$. This
implies that $b^\ell\not\in N_1$. Thus $y_0\in N_1\cap
N_1b^\ell=\varnothing$, a contradiction. Thus $M_2\times G_1$ is an irredundant member of $\mathcal{C}$. By a similar argument,
one may prove that each  $G_1\times N_j$ is an irredundant member for
the cover (for each $j=2,\dots, n$). This completes the proof of
irredundancy of the cover $\mathcal{C}$.\\
Now by Lemma 2.2(b) of \cite{BFS}, $D_1=\bigcap_{i=2}^m M_i$,
$D_2=\bigcap_{j=2}^{n}N_j$ and
$$\bigcap_{S\in\mathcal{C}}S=\big(\bigcap_{i=2}^{m} M_i\times
G_2\big)  \bigcap \big(\bigcap_{j=2}^n G_1\times N_j\big).$$ It
follows that the intersection of the cover $\mathcal{C}$ is
$D_1\times D_2$. \\
For the last statement, note that it follows from the hypothesis
that $|G_1:M_1|=|G_2:N_1|=p$, as we are assuming both $M_1$ and
$N_1$ are maximal. Thus $(M_1\times N_1)\langle (a,b)\rangle$ has
prime index $p$ in $G_1\times G_2$ and so it is a maximal
subgroup of $G_1\times G_2$. The other members of the cover
$\mathcal{C}$ are clearly maximal subgroups of $G_1\times G_2$.
This completes the proof of the last statement. $\hfill\Box$
\begin{cor}\label{c1}
  Let $V(n_i,q)$ {\rm(}$i=1,2${\rm)}
be the $n_i$-dimensional vector space over the finite field of
order $q$. If $\mathcal{C}_1=\{M_1,\dots,M_m\}$ and
$\mathcal{C}_2=\{N_1,\dots,N_n\}$ be irredundant $m$- and
$n$-covers for the abelian groups $V(n_1,q)$ and $V(n_2,q)$,
respectively; and assume that $a$ and $b$ are arbitrary elements
of $G_1\backslash M_1$ and $G_2\backslash N_1$, respectively. Then
$$\mathcal{C}=\{(M_1\times N_1)\langle (a,b)\rangle, M_i\times V(n_2,q),V(n_1,q)\times
N_j \;|\; i=2,\dots,m; \; j=2,\dots,n\}$$ is an irredundant
$(n+m-1)$-cover for $V(n_1,q)\times V(n_2,q)$ with intersection
$D_1\times D_2$, where $D_1=\cap_{i=1}^m M_i$ and
$D_2=\cap_{j=1}^n N_i$. In particular, if  all  members of
$\mathcal{C}_1$ and $\mathcal{C}_2$ are hyperplanes of the
corresponding spaces, then all members of  $\mathcal{C}$ are
hyperplanes of $V(n_1,q)\times V(n_2,q)$.
\end{cor}

\noindent{\bf Proof of Corollary \ref{cor}.}
It follows from Proposition 2.2 of \cite{AAM} and Corollary
\ref{c1}. $\hfill\Box$

\noindent{\bf Acknowledgement.}  The author is grateful to  the referee for his/her helpful comments. This research was in part supported by a grant from IPM (No. 87200118). This research was also partially
supported by  the Center of Excellence for Mathematics, University
of Isfahan.


\begin{thebibliography}{99}
\bibitem{AAM} A. Abdollahi, M.J. Ataei and A. Mohammadi Hassanabadi, Minimal blocking sets in
$\mathrm{PG}(n,2)$
 and covering groups by subgroups, to appear in  Comm.  Algebra, {\bf 36} (2008).
\bibitem{BFS} R. A. Bryce,  V. Fedri, L. Serena,  Covering group with subgroups. Bull. Austral.
Math. Soc.,   Bull. Austral. Math. Soc., {\bf 55} (1997) 469-476.
\bibitem{T} M. J. Tomkinson, Groups covered by finitely many cosets
or subgroups, Comm. Algebra, {\bf 15} (1987) No. 4, 845-859.
\end{thebibliography}
\end{document}